\documentclass[reqno,12pt]{amsart}
\usepackage{amssymb}
\usepackage{amsmath}
\usepackage{amsfonts}
\usepackage[francais]{babel}

\newtheorem{theorem}{Théorème}[section]
\theoremstyle{plain}

\newtheorem{corollary}{Corollaire}[section]

\newtheorem{problem}{Problème}[section]
\newtheorem{proposition}{Proposition}[section]
\newtheorem{remark}{Remarque}[section]

\numberwithin{equation}{section}
\input{tcilatex}

\begin{document}
\title{K(X,Y) comme sous-espace compl\'{e}ment\'{e} de $\mathcal{L}(X,Y)$}
\author{Daher Mohammad}
\maketitle

\begin{abstract}
Soient $X$,$Y$ des espaces de Banach. Dans ce travail nous montrons que $%
K(X,Y)$ contient une copie compl\'{e}ment\'{e}e de $c_{0},$ si $Y$ contient $%
\ell ^{\infty }$ isomorphiquement, ou $Y$ contient $c_{0}$ isomorphiquement
et toute suite born\'{e}e dans $Y^{\ast }$ admet une sous-suite qui converge
pr\'{e}faiblement.

Dans la suite, nous prouvons qu'il existe un espace de Banach $X$ tel que $%
K(X)$ est compl\'{e}ment\'{e} dans $\mathcal{L}(X)$, et $K(X)$ n'est pas
compl\'{e}ment\'{e} dans son bidual.

{\footnotesize Abstract. }Let $X$,$Y$ be Banach spaces. In this work, we
show that $K(X,Y)$ contains a complemented copy of $c_{0},$ if $Y$ contains
a copy of $\ell ^{\infty }$, or $Y$ contains a copy of $c_{0}$ and every
bounded sequence in $Y^{\ast }$ has a subsequence which is w$^{\ast }-$%
convergente

In the follwing, we show that there exists a Banach space $X$ such that $%
K(X) $ is complemented space in $\mathcal{L}(X)$ but $K(X)$ is not
complemented in its bidual.
\end{abstract}

Classification: 46B07, Secondaire 47B10

Mots cl\'{e}s: projection dans $\mathcal{L}(X,Y).$

\begin{center}
\textbf{Introduction.}
\end{center}

Soient $X$,$Y$ des espaces de Banach . On d\'{e}signe par $\mathcal{L}(X,Y)$
l'espace des op\'{e}rateurs born\'{e}s de $X$ dans $Y$ et par $K(X,Y)$ le
sous-espace de $\mathcal{L}(X,Y)$ form\'{e} des op\'{e}rateurs compacts$.$

Dans ce travail, nous montrons que $K(X,Y)$ contient une copie compl\'{e}ment%
\'{e}e de $c_{0},$ si $Y$ contient $\ell ^{\infty }$ isomorphiquement, ou $Y$
contient $c_{0}$ isomorphiquement et toute suite born\'{e}e dans $Y^{\ast }$
admet une sous-suite qui converge pr\'{e}faiblement. Dans la suite, nous
retrouvons des r\'{e}sultats de M.Feder \cite{Fe}, J.Johnson \cite{Joh}, et
G.Emmanuele \cite{Emm} concernant l'existence d'une projection continue de $%
\mathcal{L}(X,Y)$ sur $K(X,Y).$

Finalement, nous racact\'{e}risons l'existence d'une projection continue:$%
K(X,Y)\rightarrow K(X,Y)^{\ast \ast },$ si $Y$ ou $X^{\ast }$ a la propri%
\'{e}t\'{e} de l'approximation born\'{e}e.

Pour tout Banach $X,$ on d\'{e}signe par $B_{X}$ la boule unit\'{e} ferm\'{e}
de $X$ et par $X^{\ast }$ le dual de $X.$

Notons pour tout $x\in X$ et tout $x^{\ast }\in X^{\ast }$ $(x,x^{\ast
})=x^{\ast }(x).$

D\'{e}signons par $(e_{k})_{k\geq 0}$ la base canonique de $c_{0}.$

Dans la suite on fixe deux espaces de Banach $X,Y$ de dimensions infinies.

\begin{proposition}
\label{kj}Supposons que $Y$ contient $c_{0}$ isomorphiquement. Alors il
existe un op\'{e}rateur $\sigma :\ell ^{\infty }\rightarrow \mathcal{L}(X,Y)$
tel que
\end{proposition}

a) $\ell ^{\infty }\approx \sigma (\ell ^{\infty }).$

b) $\sigma (\ell ^{\infty })\cap K(X,Y)=\sigma (c_{0}).$

c) Il existe une projection continue $P:\mathcal{L}(X,Y)\rightarrow \sigma
(\ell ^{\infty })$ telle que $P_{\mid _{K(X,Y_{1})}}$ est une projection sur 
$\sigma (c_{0}),$ o\`{u} $Y_{1}$ est un sous-espace de Banach de $Y$
isomorphe \`{a} $c_{0}.$

d) Supposons $Y$ contient $\ell ^{\infty }$ isomorphiquement, ou toute suite
born\'{e}e dans $Y^{\ast }$ admet une sous-suite qui converge pr\'{e}%
faiblement. Alors il existe $\Lambda \subset \mathbb{N}$ et une projection
continue $Q:\mathcal{L}(X,Y)\rightarrow \sigma (\ell ^{\infty }(\Lambda ))$
telle que $Q_{\mid _{K(X,Y)}}$ est une projection sur $\sigma (c_{0}(\Lambda
)).$

D\'{e}monstration.

a).

Soient $Y_{1}$ un sous-espace de Banach de $Y$ et $U:c_{0}\rightarrow Y_{1}$
un tel isomorphisme. D'apr\`{e}s \cite[Chap.XII]{Di}-\cite{Jos}-\cite{Ness},
il existe une suite $(x_{k}^{\ast })_{k\geq 0}$ dans la sph\`{e}re unit\'{e}
de $X^{\ast }$ telle que $x_{k}^{\ast }\underset{k\rightarrow +\infty }{%
\rightarrow }0$ pr\'{e}faiblement. On d\'{e}finit l'op\'{e}rateur $\sigma
:\ell ^{\infty }\rightarrow \mathcal{L}(X,Y)$ par $\sigma (\alpha )(x)=%
\underset{k\geq 0}{\dsum }\alpha _{k}(x,x_{k}^{\ast })U(e_{k}),$ $\alpha
=(\alpha _{k})_{k\geq 0}\in \ell ^{\infty },$ $x\in X.$

Pour tout $k\in \mathbb{N}$, il existe $x_{k}\in X$ tel que $%
(x_{k},x_{k}^{\ast })\neq 0.$ Cela implique que $\sigma $ est injectif$.$

\emph{Etape 1:}\textbf{\ }Montrons que $\sigma $ est un op\'{e}rateur born%
\'{e}.

Pour tout $\alpha \in \ell ^{\infty }$ et tout $x\in X$ on a

\begin{eqnarray*}
\left\Vert \sigma (\alpha )(x)\right\Vert _{Y} &=&\left\Vert \underset{k\geq
0}{\dsum }\alpha _{k}(x,x_{k}^{\ast })U(e_{k})\right\Vert _{Y}=\left\Vert 
\underset{}{U(\dsum_{k\geq 0}}\alpha _{k}(x,x_{k}^{\ast })e_{k})\right\Vert
_{Y}\leq \\
\left\Vert U\right\Vert \left\Vert \underset{k\geq 0}{\dsum }\alpha
_{k}(x,x_{k}^{\ast })e_{k}\right\Vert _{c_{0}} &\leq &\left\Vert
U\right\Vert \left\Vert \alpha \right\Vert _{\ell ^{\infty }}\times
\sup_{k\geq 0}\left\vert (x,x_{k}^{\ast })\right\vert \leq \left\Vert
U\right\Vert \times \left\Vert \alpha \right\Vert _{\ell ^{\infty }}\times
\left\Vert x\right\Vert .
\end{eqnarray*}

Il en r\'{e}sulte que $\sigma $ est un op\'{e}rateur born\'{e}.

\emph{Etape 2:}\textbf{\ }Montrons que $\sigma $ est un isomorphisme sur son
image.

Fixons $k\in \mathbb{N}.$ D\'{e}signons par $z_{k}^{\ast }$ la forme lin\'{e}%
aire d\'{e}finie sur $Y_{1}$ par $(Ue_{j},z_{k}^{\ast })=\delta _{jk},$ $%
j\in \mathbb{N}$, o\`{u} $\delta _{jk}$ est le symbole de Kronecker$.$

\bigskip

Soient $\alpha \in \ell ^{\infty }$, $\varepsilon >0$ et $k\in \left\{ j\in 
\mathbb{N}^{{}};\text{ }\alpha _{j}\neq 0\right\} $. Il existe $x_{k}\in
B_{X}$ tel que $1=\left\Vert x_{k}^{\ast }\right\Vert \leq \left\vert
(x_{k},x_{k}^{\ast })\right\vert +\varepsilon /\left\vert \alpha
_{k}\right\vert .$ Pour tout $k\in \left\{ j\in \mathbb{N}\text{; }\alpha
_{j}\neq 0\right\} $ on a%
\begin{eqnarray*}
\left\vert \alpha _{k}\right\vert &\leq &\left\vert \alpha _{k}\right\vert
\times \left[ \left\vert (x_{k},x_{k}^{\ast })\right\vert +\varepsilon
/\left\vert \alpha _{k}\right\vert \right] =\left\vert \alpha
_{k}\right\vert \times \left\vert (x_{k},x_{k}^{\ast })\right\vert
+\varepsilon \\
&=&\left\vert (\sigma (\alpha )(x_{k}),z_{k}^{\ast })\right\vert
+\varepsilon \leq \left\Vert \sigma (\alpha )\right\Vert \times \sup_{k\geq
0}\left\Vert z_{k}^{\ast }\right\Vert _{Y_{1}^{\ast }}+\varepsilon .
\end{eqnarray*}

$.$ Cela entra\^{\i}ne que $\sigma $ est un isomorphisme sur son image.

b).

\bigskip Montrons que $\sigma (c_{0})\subset \sigma (\ell ^{\infty })\cap
K(X,Y).$

Soient $\alpha \in c_{0}$ et $\varepsilon >0.$ Il existe $k_{0}\in \mathbb{N}
$ tel que $\left\vert \alpha _{k}\right\vert \leq \varepsilon $ pour tout $%
k\geq k_{0}$. Consid\'{e}rons $T_{n}:X\rightarrow Y$ \ l'op\'{e}rateur d\'{e}%
finie par%
\begin{equation*}
T_{n}(x)=\overset{n}{\underset{k=0}{\dsum }}\alpha _{k}(x,x_{k}^{\ast
})Ue_{k},\text{ }x\in X,\text{ }n\in \mathbb{N}.
\end{equation*}

Il est \'{e}vident que la suite $(T_{n})_{n\geq 0}$ est dans $K(X,Y)$ et que 
$\left\Vert T_{n}-\sigma (\alpha )\right\Vert _{\mathcal{L}%
(X,Y)}<\varepsilon $ pour tout $n\geq k_{0},$ donc $\sigma (\alpha )$ est un
op\'{e}rateur compact.

Montrons que $\sigma (\ell ^{\infty })\cap K(X,Y)\subset \sigma (c_{0}).$

Remarquons d'abord que si $\alpha _{1},...,\alpha _{n}\in \mathbb{C}$ et $y=%
\underset{k\leq n}{\dsum }\alpha _{k}U(e_{k}),$ $(y,z_{j}^{\ast })\underset{%
j\rightarrow \infty }{\rightarrow }0,$ donc par densit\'{e} pour tout $y\in
Y_{1}$ $(y,z_{j}^{\ast })\underset{j\rightarrow \infty }{\rightarrow }0.$

Soit $\alpha \in \ell ^{\infty }$ tel que $\sigma (\alpha )\in K(X,Y).$
Supposons que $\alpha \notin c_{0}$. Il existe donc $\varepsilon _{0}>0,$
tel que pour tout $k\in \mathbb{N},$ il existe $m_{k}\in \mathbb{N}$ v\'{e}%
rifiant 
\begin{equation}
m_{k}\geq k\text{ et }\left\vert \alpha _{m_{k}}\right\vert \geq \varepsilon
_{0}.  \label{c}
\end{equation}

Posons $T=\sigma (\alpha )\in K(X,Y_{1}).$ Pour tout $k\in \mathbb{N},$ il
existe $x_{k}$ dans la boule unit\'{e} de $X$ tel que $\left\vert
(x_{k},x_{k}^{\ast })\right\vert \geq 1-1/k+2.$ Il est clair que $\alpha
_{k}=(Tx_{k},z_{k}^{\ast })/(x_{k},x_{k}^{\ast })$ pour tout $k\in \mathbb{N}%
.$

Comme $T$ est un op\'{e}rateur compact, il existe une sous-suite $%
(x_{m_{k_{j}}})_{j\geq 0}$ dans $X$ telle que $T(x_{m_{k_{j}}})\underset{%
j\rightarrow +\infty }{\rightarrow }y\in Y_{1}.$ D'autre part, pour tout $%
j\in \mathbb{N}$ 
\begin{equation*}
\alpha _{m_{k_{j}}}=\left[ (Tx_{m_{k_{j}}}-y,z_{m_{k_{j}}}^{\ast
})+(y,z_{m_{k_{j}}}^{\ast })\right] /(x_{m_{k_{j}}},x_{m_{k_{j}}}^{\ast })
\end{equation*}

et $(y,z_{m_{k_{j}}}^{\ast })\underset{j\rightarrow +\infty }{\rightarrow }%
0. $ Par cons\'{e}quent $\alpha _{m_{k_{j}}}\underset{j\rightarrow +\infty }{%
\rightarrow }0$, car $\left\vert (x_{k},x_{k}^{\ast })\right\vert \geq
1-1/k+2$ pour tout $k\in \mathbb{N},$ ce qui est impossible d'apr\`{e}s (\ref%
{c})$.$ Donc $\alpha \in c_{0}.$

\bigskip c).

D'apr\`{e}s le th\'{e}or\`{e}me de Hahn-Banach, pour tout $k\in \mathbb{N},$
il existe $y_{k}^{\ast }\in Y^{\ast }$ qui prolonge $z_{k}^{\ast }\in
Y_{1}^{\ast }$ et $\left\Vert z_{k}^{\ast }\right\Vert _{Y_{1}^{\ast
}}=\left\Vert y_{k}^{\ast }\right\Vert _{Y^{\ast }}$

\bigskip Soit maintenant $T\in \mathcal{L}(X,Y).$ D\'{e}finissons $a(T)\in
\ell ^{\infty }$ par $\alpha _{k}(T)=(Tx_{k},y_{k}^{\ast
})/(x_{k},x_{k}^{\ast }),$ $k\in \mathbb{N}$ et $P:\mathcal{L}%
(X,Y)\rightarrow \sigma (\ell ^{\infty })$ par $P(T)=\sigma \left[ \alpha (T)%
\right] ,$ $T\in \mathcal{L}(X,Y).$

\emph{Etape 1:} Montrons que $P$ est une projection.

Soient $\alpha \in \ell ^{\infty }$. Posons $T=\sigma (\alpha ).$ Remarquons
que%
\begin{eqnarray*}
\alpha _{k}(T) &=&(Tx_{k},y_{k}^{\ast })/(x_{k},x_{k}^{\ast })= \\
(Tx_{k},z_{k}^{\ast })/(x_{k},x_{k}^{\ast }) &=&\alpha
_{k}(U(e_{k}),z_{k}^{\ast })(x_{k},x_{k}^{\ast })/(x_{k},x_{k}^{\ast
})=\alpha _{k}.
\end{eqnarray*}

Donc $P$ est une projection.

\emph{Etape 2:} Montrons que $P$ est continue.

Pour tout $T\in \mathcal{L}(X,Y)$ $\ $nous avons $\left\Vert \sigma (\alpha
(T))\right\Vert \leq \left\Vert \sigma \right\Vert \times \left\Vert \alpha
(T)\right\Vert _{\ell ^{\infty }}$ $\leq \left\Vert \sigma \right\Vert
\times \sup_{k\geq 0}\left\vert \left[ (Tx_{k},y_{k}^{\ast
})/(x_{k},x_{k}^{\ast })\right] \right\vert \leq \left\Vert \sigma
\right\Vert \times \left[ \left\Vert T\right\Vert \right] \times \left[
\sup_{k\geq 0}(\left\Vert y_{k}^{\ast }\right\Vert /(1-1/k+2))\right]
.\blacksquare $

\emph{Etape 3: }Soit $T\in K(X,Y_{1})$. Montrons que $\alpha (T)\in c_{0}.$

Comme $P$ est une projection, d'apr\`{e}s b), $T=\sigma (\alpha (T))\in
K(X,Y)\cap \sigma (\ell ^{\infty })=\sigma (c_{0}),$ c'est-\`{a}-dire que $%
\alpha (T)\in c_{0}.\blacksquare $

d).

\emph{Cas 1: }$Y$ contient $\ell ^{\infty }$ isomorphiquement.

Comme $\ell ^{\infty }$ est injectif, on peut supposer que $Y\approx \ell
^{\infty }.$ Consid\'{e}rons $G:Y\rightarrow \ell ^{\infty }$ un
isomorphisme. Notons $(f_{n})_{n\geq 0}$ la base canonique de $\ell ^{1},$ $%
h_{n}^{\ast }=G^{\ast }f_{n}\in Y^{\ast },$ $n\in \mathbb{N}$ et $%
U:c_{0}\rightarrow Y_{1},$ la restriction de $G^{-1}$ \`{a} $c_{0}.$

Soit maintenant $T\in \mathcal{L}(X,Y).$ D\'{e}finissons $\beta (T)\in \ell
^{\infty }$ par $\beta _{k}(T)=(Tx_{k},h_{k}^{\ast })/(x_{k},x_{k}^{\ast }),$
$k\in \mathbb{N}$ et $Q:\mathcal{L}(X,Y)\rightarrow \sigma (\ell ^{\infty })$
par $Q(T)=\sigma \left[ \beta (T)\right] ,$ $T\in \mathcal{L}(X,Y).$ Par un
argument analogue \`{a} celui de c), on montre que $Q$ est une projection
continue et $Q_{\mid _{K(X,Y)}}$ est une projection sur $\sigma (c_{0})$
(ici $\Lambda =\mathbb{N}).\blacksquare $

\emph{Cas 2: }Toute suite born\'{e}e dans $Y^{\ast }$ admet une sous-suite pr%
\'{e}faiblement convergeante.

Il existe une sous-suite $(y_{n_{k}}^{\ast })_{k\geq 0}$ telle que $%
y_{n_{k}}^{\ast }\underset{k\rightarrow \infty }{\rightarrow }y^{\ast }$ pr%
\'{e}faiblement dans $Y^{\ast }$. D\'{e}singons par $\Lambda =\left\{ n_{k};%
\text{ }k\in \mathbb{N}\right\} .$

Soit $U:c_{0}(\Lambda )\rightarrow Y_{1}$ un isomorphisme. On d\'{e}finit
l'op\'{e}rateur $\sigma :\ell ^{\infty }(\Lambda )\rightarrow \mathcal{L}%
(X,Y)$ par $\sigma (\alpha )(x)=\underset{k\in \Lambda }{\dsum }\alpha
_{k}(x,x_{k}^{\ast })U(e_{k}),$ $x\in X,\alpha \in \ell ^{\infty }(\Lambda )$

Pour tout $T\in \mathcal{L}(X,Y)$ on d\'{e}finit $\gamma (T)\in \ell
^{\infty }(\Lambda )$ par $\gamma _{k}(T)=(Tx_{k},y_{k}^{\ast }-y^{\ast
})/(x_{k},x_{k}^{\ast }),$ $k\in \Lambda $ et $Q:\mathcal{L}(X,Y)\rightarrow
\sigma (\ell ^{\infty }(\Lambda ))$ par $Q(T)=\sigma \left[ \gamma (T)\right]
,$ $T\in \mathcal{L}(X,Y).$

Montrons que $Q$ est une projection.

Consid\'{e}rons $T=\sigma (\alpha ),$ $\alpha \in \ell ^{\infty }(\Lambda )$%
. Pour tout $k\in \Lambda $ on a $\gamma _{k}(T)=(Tx_{k},y_{k}^{\ast
}-y^{\ast })/(x_{k},x_{k}^{\ast })=(Tx_{k},z_{k}^{\ast }-y^{\ast
})/(x_{k},x_{k}^{\ast })=\left[ \underset{j\in \Lambda }{\dsum }\alpha
_{j}(x_{k},x_{j}^{\ast })(U(e_{j}),z_{k}^{\ast }-y^{\ast })\right]
/(x_{k},x_{k}^{\ast })=\alpha _{k},$ car pour tout $j\in \Lambda $ $%
(U(e_{j}),y^{\ast })=\lim_{k\rightarrow \infty }(U(e_{j}),z_{n_{k}}^{\ast
})=0.$

Par un argument analogue \`{a} celui de c), on montre que $Q$ est continue.

Montrons finalement que $\gamma (T)\in c_{0}(\Lambda ),$ pour tout $T\in
K(X,Y).$

Soit $T\in K(X,Y).$ Supposons que $\gamma (T)\notin c_{0}$, il existe $%
\varepsilon _{0}>0$ tel que pour tout $k\in \mathbb{N},$ il existe $%
m_{k}^{{}}\in \Lambda $ v\'{e}rifiant 
\begin{equation}
m_{k}^{{}}\geq k\text{ et }\left\vert \gamma _{m_{k}^{{}}}(T)\right\vert
\geq \varepsilon _{0}.  \label{iv}
\end{equation}

D'autre part, il existe une sous-suite $(x_{m_{k_{j}}})_{j\geq 0}$ telle que 
$Tx_{m_{k_{j}}}\underset{j\rightarrow \infty }{\rightarrow }y\in Y.$

Pour tout $j\in \mathbb{N}$ on a

\begin{equation*}
\gamma _{m_{k_{j}}}(T)=\left[ (Tx_{m_{k_{j}}}-y,y_{m_{k_{j}}}^{\ast
}-y^{\ast })+(y,y_{m_{k_{j}}}^{\ast }-y^{\ast })\right]
/(x_{m_{k_{j}}},x_{m_{k_{j}}}^{\ast }).
\end{equation*}

Il est clair que $\gamma _{m_{k_{j}}}(T)\underset{j\rightarrow \infty }{%
\rightarrow }0,$ car $(y,y_{m_{k_{j}}}-y^{\ast })\underset{j\rightarrow
\infty }{\rightarrow }0,$ ceci est impossible, donc $\gamma (T)\in
c_{0}(\Lambda ).\blacksquare $

Les deux corollaires suivants sont une cons\'{e}quence de la proposition \ref%
{kj}-d.

\begin{corollary}
\label{nu}Supposons $Y$ contient $\ell ^{\infty }$ isomorphiquement. Alors $%
K(X,Y)$ contient une copie compl\'{e}ment\'{e}e de $c_{0}.$
\end{corollary}

\begin{corollary}
\label{mu}Supposons que $Y$ contient $c_{0}$ isomorphiquement et que toute
suite born\'{e}e dans $Y^{\ast }$ admet une sous-suite qui converge pr\'{e}%
faiblement. Alors $K(X,Y)$ contient une copie compl\'{e}ment\'{e}e de $%
c_{0}. $
\end{corollary}

\begin{corollary}
\label{jj}\cite[Coroll.4]{Fe}-\cite[Coroll.4]{Ch}-\cite[Th. 4]{Joh}Supposons
que $Y$ contient $c_{0}$ isomorphiquement. Alors il n'existe aucune
projection continue:$\mathcal{L(}X,Y)\rightarrow K(X,Y)$.
\end{corollary}

D\'{e}monstration$.$

Suppososns qu'il existe une projection continue $P:\mathcal{L}%
(X,Y)\rightarrow K(X,Y).$

\emph{Cas 1}\textbf{:} $Y$ ne contient pas $\ell ^{\infty }$
isomorphiquement et $X$ ne contient pas une copie compl\'{e}ment\'{e}e de $%
\ell ^{1}.$

D'apr\`{e}s \cite[Th.4]{Kalt}, $K$($X,Y)$ ne contient pas $\ell ^{\infty }$
isomorphiquement$.$ D'autre part, la proposition \ref{kj}, nous montre qu'il
existe un isomorphisme $\sigma :\ell ^{\infty }\rightarrow \sigma (\ell
^{\infty })\subset \mathcal{L}(X,Y)$. En appliquant \cite[Prop.2]{Kalt}, on
voit que $P\circ \sigma :\ell ^{\infty }\rightarrow K(X,Y)$ est faiblement
compact, ceci est impossible, car $P\circ \sigma _{\mid _{c_{0}}}=\sigma
_{\mid _{c_{0}}}$un isomorphisme.$\blacksquare $

\emph{Cas 2}\textbf{: }$Y$ ne contient pas $\ell ^{\infty }$
isomorphiquement $\ $et $X$ contient une copie compl\'{e}ment\'{e}e de $\ell
^{1}.$

Ce cas implique qu'il existe une porojection continue $P_{1}:\mathcal{L}%
(\ell ^{1},Y)\rightarrow K(\ell ^{1},Y).$ D'apr\`{e}s \cite[Th.6]{Kalt}, $%
K(\ell ^{1},Y)$ ne contient pas $c_{0}$ isomorphiquement$,$ ce qui est
impossible, car $K(\ell ^{1},Y)$ contient $\ell ^{\infty }$ isomorphiquement$%
.\blacksquare $

\emph{Cas 3}\textbf{: }$Y$ contient $\ell ^{\infty }$ isomorphiquement$.$
D'apr\`{e}s la proposition \ref{kj}-d), il existe une projection $%
Q:K(X,Y)\rightarrow \sigma (c_{0}).$ Comme $\ell ^{\infty }$ est un espace
de Grothendieck, d'apr\`{e}s \cite[Groth]{Groth} $\sigma \circ Q\circ P:\ell
^{\infty }\rightarrow \sigma (c_{0})$ est faiblement compact, ce qui est
impossible, car $\sigma \circ Q\circ P_{\mid _{c_{0}}}$un isomorphisme.$%
\blacksquare $

\begin{corollary}
\label{dfg}\cite[Th.2]{Emm}-\cite[Coroll.7]{Ch}-\cite[Croll.7]{Joh}Supposons
qu'il existe une projection continue $P:\mathcal{L}(X,Y)\rightarrow K(X,Y).$
Alors $K(X,Y)$ ne contient pas $c_{0}$ isomorphiquement$.$
\end{corollary}

D\'{e}monstration.

\bigskip D'apr\`{e}s le corollaire \ref{jj}, $Y$ ne contient pas $c_{0}$
isomorphiquement.

Montrons d'abord que $K(X,Y)$ ne contient pas $\ell ^{\infty }$
isomorphiquement.

Supposons que $K(X,Y)$ contient $\ell ^{\infty }$ isomorphiquement$.$ Le r%
\'{e}sultat de \cite[Th.4]{Kalt}, nous indique que $X$ contient une copie
compl\'{e}ment\'{e}e de $\ell ^{1},$ il en r\'{e}sulte qu'il existe donc une
projection continue:$\mathcal{L}(\ell ^{1},Y)\rightarrow K(\ell ^{1},Y),$ ce
qui est impossible d'apr\`{e}s \cite[Th.6]{Kalt}. On en d\'{e}duit que $%
K(X,Y)$ ne contient pas $\ell ^{\infty }$ isomorphiquement$.$

Supposons qu'il existe un isomorphisme $\sigma _{1}:c_{0}\rightarrow \sigma
_{1}(c_{0})\subset K(X,Y).$ Pour tout $\alpha \in \ell ^{\infty }$ et tout $%
x\in X$ on a%
\begin{eqnarray*}
sup\left\{ \left\Vert \underset{J}{\dsum }\alpha _{n}\sigma
_{1}(e_{n})x\right\Vert _{Y};\text{ }J\text{ \ est un sous-ensemble de }%
\mathbb{N}\text{ fini}\right\}  &<& \\
\left\Vert x\right\Vert \times sup\left\{ \left\Vert \underset{J}{\dsum }%
\alpha _{n}\sigma _{1}(e_{n})\right\Vert _{\mathcal{L}(X,Y)};\text{ }J\text{
\ est un sous-ensemble de }\mathbb{N}\text{ fini}\right\}  &<& \\
\left\Vert x\right\Vert \times \left\Vert \sigma _{1}\right\Vert sup\left\{
\left\Vert \underset{J}{\dsum }\alpha _{n}e_{n}\right\Vert _{c_{0}};\text{ }J%
\text{ \ est un sous-ensemble de }\mathbb{N}\text{ fini}\right\}  &<&+\infty
.
\end{eqnarray*}

Comme $Y$ ne contient pas $c_{0}$ isomorphiquement$,$ la s\'{e}rie $\underset%
{n\geq 0}{\dsum }\alpha _{n}\sigma _{1}(e_{n})x$ converge dans $Y$ pour tout 
$x\in X$ \cite[prop.3]{Kalt}$.$ On d\'{e}finit $\sigma _{2}:\ell ^{\infty
}\rightarrow \mathcal{L}(X,Y)$ par $\sigma _{2}(\alpha )(x)=\underset{k\geq 0%
}{\dsum }\alpha _{k}\sigma _{1}(e_{k})(x),$ $x\in X$ $.$ Consid\'{e}rons $%
P\circ \sigma _{2}:\ell ^{\infty }\rightarrow K(X,Y)$. D'apr\`{e}s \cite[%
Prop2]{Kalt}, $P\circ \sigma _{2}$ est faiblement compact, car $K(X,Y)$ ne
contient pas $\ell ^{\infty }$ isomorphiquement$,$ c'est impossible, car $%
P\circ \sigma _{2_{\mid _{c_{0}}}}=\sigma _{1}.\blacksquare $

\begin{proposition}
\label{nsf}Supposons que $\mathcal{L}(X,Y)$ ne contient pas $\ell ^{\infty }$
isomrphiquement. Alors $K(X,Y)$ ne contient pas $c_{0}$ isomorphiquement.
\end{proposition}

D\'{e}monstration.

Supposons qu'il existe un isomorphisme $\sigma _{1}:c_{0}\rightarrow \sigma
_{1}(c_{0})\subset K(X,Y).$ Remarquons d'apr\`{e}s la proposition \ref{kj}
que $Y$ ne contient pas $c_{0}$ isomorphiquement et que 
\begin{equation*}
sup\left\{ \left\Vert \underset{}{\underset{J}{\dsum }\alpha _{n}\sigma
_{1}(e_{n})}(x)\right\Vert _{Y};\text{ }J\text{ est un sous-ensemble fini de 
}\mathbb{N}\right\} <+\infty 
\end{equation*}%
pour tout $\alpha \in \ell ^{\infty }$ et tout $x\in X.$ L'espace $Y$ ne
contient pas $c_{0}$ isomorphiquement, donc d'apr\`{e}s \cite[Prop.3]{Kalt}, 
$\underset{}{\text{la s\'{e}rie }\underset{k\geq 0}{\dsum }}\alpha
_{n}\sigma _{1}(e_{n})x$ converge dans $Y$ pour tout $x\in X.$ On d\'{e}%
finit $\sigma _{2}:\ell ^{\infty }\rightarrow \mathcal{L}(X,Y)$ par $\sigma
_{2}(\alpha )x=\underset{n\geq 0}{\dsum }\alpha _{n}\sigma _{1}(e_{n})x,$ $%
\alpha \in \ell ^{\infty }$, $x\in X.$ Comme $\mathcal{L}(X,Y)$ ne contient
pas $\ell ^{\infty }$ isomorphiquement, d'apr\`{e}s \cite[Prop.2]{Kalt} $%
\sigma _{2}$ est faiblement compact. Il en r\'{e}sulte que $\sigma _{1}$ est
faiblement compact, ceci est impossible, par cons\'{e}quent $K(X,Y)$ ne
contient pas $c_{0}$ isomorphiquement.$\blacksquare $

\bigskip

Soit $X,Y$ deux espaces de Banach. D\'{e}signons par $X\overset{\vee }{%
\otimes }Y$ (resp. $X\widehat{\otimes }Y)$ le produit injectif de $X$,$Y$
(resp. le produit projectif de $X,Y).$ D\'{e}signons d'autre part, par $%
B(X\times Y)$ l'espace des formes bilin\'{e}aires $u:X\times Y\rightarrow 
\mathbb{C}$ telles que sup$\left\{ \left\vert u(x,y)\right\vert ;\text{ }%
\left\Vert x\right\Vert \text{,}\left\Vert y\right\Vert \leq 1\right\}
<\infty $.

\begin{theorem}
\label{ert}Il existe un espace de Banach $Z$ qui poss\`{e}des des propri\'{e}%
t\'{e}s suivantes:
\end{theorem}

1) $K(Z)$ est un espace compl\'{e}ment\'{e} de $\mathcal{L}(Z).$

2) $K(Z)$ n'est pas compl\'{e}ment\'{e} dans son bidual.

\begin{proposition}
\label{rutt}Supposons que $Y$ ou $X^{\ast }$ a la propri\'{e}t\'{e} de
l'approximation born\'{e}e. Alors $K(X,Y)$ est complement\'{e} dans $%
K(X,Y)^{\ast \ast }$, si et seulement si $K(X,Y)$ et $Y$ sont compl\'{e}ment%
\'{e}s dans $\mathcal{L}(X,Y)$ et dans $Y^{\ast \ast }$ respectivement.
\end{proposition}

D\'{e}monstration de proposition \ref{rutt}.

\bigskip Supposons qu'il existe une projection continue $Q:K(X,Y)^{\ast \ast
}\rightarrow K(X,Y)$ $.$ Montrons que $Y$ est compl\'{e}ment\'{e} dans son
bidual.

Consid\'{e}rons $x_{0}$ dans la sph\`{e}re unit\'{e} de $X$ et $x^{\ast }$
dans la sph\`{e}re unit\'{e} de $X^{\ast }$ tels que $(x_{0},x^{\ast })=1.$
On d\'{e}finit l'op\'{e}rateur $H:Y\rightarrow K(X,Y)$ par $%
H(y)(x)=(x,x^{\ast })y,$ $x\in X,y\in Y$ et l'op\'{e}rateur $%
G:K(X,Y)\rightarrow Y,$ par $G(T)=Tx_{0},$ $T\in K(X,Y).$ Notons $%
Q_{1}=G\circ Q\circ H^{\ast \ast }:Y^{\ast \ast }\rightarrow Y.$ Il est
clair que $Q_{1}$ est une projection continue sur $Y.$

Montrons que $K(X,Y)$ est compl\'{e}ment\'{e} dans $\mathcal{L}(X,Y).$

\emph{Cas 1: }$Y$ a la propri\'{e}t\'{e} de l'approximation born\'{e}e.

D'apr\`{e}s \cite[Lemma 2]{Joh}, il existe un isomorphisme $U:\mathcal{L}%
(X,Y)\rightarrow K(X,Y)^{\ast \ast }$ tel que la restriction de $U$ \`{a} $%
K(X,Y)$ est l'identit\'{e}.

Notons $P=Q\circ U:\mathcal{L}(X,Y)\rightarrow K(X,Y)$. Il est facile de
voir que $P$ est une projection continue sur $K(X,Y)$.

\emph{Cas 2:} $X^{\ast }$ a la propri\'{e}t\'{e} de l'approximation born\'{e}%
e.

Il existe une suite g\'{e}n\'{e}ralis\'{e}e $(V_{i})_{i\in I}$ d'op\'{e}%
rateurs de rang finis :$X^{\ast }\rightarrow X^{\ast }$ telle que $%
V_{i}\rightarrow I_{X^{\ast }}$ uniform\'{e}ment sur tout compact de $%
X^{\ast }$ et sup$_{i\in I}\left\Vert V_{i}\right\Vert _{\mathcal{L}(X^{\ast
})}<\infty .$ Posons pour $i\in I$ et $T\in \mathcal{L}(X,Y),$ $%
U_{i}^{T}=V_{i}\circ T^{\ast }:Y^{\ast }\rightarrow X^{\ast }$ et $%
H_{i}^{T}=Q_{1}\circ (U_{i}^{T})_{\mid _{X}}^{\ast }:$ $X\rightarrow Y.$ Il
est \'{e}vident que $H_{i}^{T}\in K(X,Y).$ On d\'{e}finit l'op\'{e}rateur $U:%
\mathcal{L}(X,Y)\rightarrow K(X,Y)^{\ast \ast }$, par $U(T)=\lim_{\mathcal{U}%
}H_{i}^{T},$ $T\in \mathcal{L}(X,Y)$ (limite pr\'{e}faible suivant $\mathcal{%
U}$ dans $K(X,Y)^{\ast \ast }).$ On se propose de montrer que $P=Q\circ U:$ $%
\mathcal{L}(X,Y)\rightarrow K(X,Y)$ est une projection continue.

Remarquons d'abord que $P$ est continue, car la suite $(V_{i})_{i\in I}$ est
born\'{e}e dans $\mathcal{L}(X^{\ast }).$ Consid\'{e}rons $T\in K(X,Y).$
Comme $T^{\ast }$ est compact, $\left\Vert U_{i}^{T}-T^{\ast }\right\Vert _{%
\mathcal{L}(Y^{\ast },X^{\ast })}\rightarrow 0,$ ceci implique que $%
\left\Vert (U_{i}^{T})^{\ast }-T^{\ast \ast }\right\Vert _{\mathcal{L}%
(X^{\ast \ast },Y^{\ast \ast })}\rightarrow 0$. Il en r\'{e}sulte que $%
\left\Vert H_{i}^{T}-Q_{1}\circ T\right\Vert _{\mathcal{L}(X,Y)}=\left\Vert
H_{i}^{T}-T\right\Vert _{\mathcal{L}(X,Y)}\rightarrow 0,$ par cons\'{e}quent 
$P(T)=Q\circ U(T)=Q(T)=T.\blacksquare $

\bigskip Supposons qu'il existe une projection $P:\mathcal{L}%
(X,Y)\rightarrow K(X,Y)$ et $Q_{1}:Y^{\ast \ast }\rightarrow Y$ une
projection continue. Montrons que $K(X,Y)$ est compl\'{e}ment\'{e} dans $%
K(X,Y)^{\ast \ast }.$

Comme $Y$ ou $X^{\ast }$ \`{a} la propri\'{e}t\'{e} de l'approximation born%
\'{e}e, $K(X,Y)=X^{\ast }\overset{\vee }{\otimes }Y.$ On d\'{e}finit l'op%
\'{e}rateur $J:X\widehat{\otimes }Y^{\ast }\rightarrow K(X,Y)^{\ast }$

par $(J\underset{k\leq n}{\dsum }x_{k}\otimes y_{k}^{\ast },\underset{j\leq m%
}{\dsum }x_{j}^{\ast }\otimes y_{j})=\underset{k\leq n,j\leq m}{\dsum }%
(x_{k},x_{j}^{\ast })\times (y_{j},y_{k}^{\ast }),$ les $x_{k}\in X,$ les $%
y_{k}^{\ast }\in Y^{\ast },$ les $x_{j}^{\ast }\in X^{\ast }$ et les $%
y_{j}\in Y.$

Montrons que $J$ est un op\'{e}rateur born\'{e}.

Consid\'{e}rons $(x_{k})_{k\leq n}\in X,$ ($y_{k}^{\ast })_{k\leq n}\in
Y^{\ast },$ $(x_{j}^{\ast })_{j\leq m}\in X^{\ast },$ $(y_{j})_{j\leq m}\in
Y.$ On a alors

\begin{eqnarray*}
&&\left\vert (J\underset{k\leq n}{\dsum }x_{k}\otimes y_{k}^{\ast },\underset%
{j\leq m}{\dsum }x_{j}^{\ast }\otimes y_{j})\right\vert \\
&=&\left\vert \underset{k\leq n}{\dsum (y_{k}^{\ast },}\left[ \underset{}{%
\text{\ }\underset{j\leq m}{\dsum }}(x_{j}^{\ast }\otimes y_{j})\right]
(x_{k}))\right\vert \\
&\leq &\underset{k\leq n}{\dsum }\left\Vert y_{k}^{\ast }\right\Vert
_{Y^{\ast }}\left\Vert \left[ \underset{j\leq m}{\dsum }x_{j}^{\ast }\otimes
y_{j})\right] (x_{k})\right\Vert _{Y} \\
&\leq &\underset{k\leq n}{\dsum }\left\Vert y_{k}^{\ast }\right\Vert
_{Y^{\ast }}\left\Vert x_{k}\right\Vert _{X}\left\Vert \underset{j\leq m}{%
\dsum }x_{j}^{\ast }\otimes y_{j})\right\Vert _{K(X,Y)}.
\end{eqnarray*}

Donc $J$ est continue.

Soit maintenant $R\in K(X,Y)^{\ast \ast }.$ On d\'{e}finit $\xi _{R}\in (X%
\widehat{\otimes }Y^{\ast })^{\ast }$ par $\xi _{R}(u)=(J(u),R)$, $u\in X%
\widehat{\otimes }Y^{\ast }.$ D'apr\`{e}s \cite[coroll.2,chap.VIII-2]{Du}, $%
\xi _{R}\in B(X\times Y^{\ast })=\mathcal{L}(X,Y^{\ast \ast }).$ Consid\'{e}%
rons $V:\mathcal{L}(X,Y^{\ast \ast })\rightarrow \mathcal{L}(X,Y)$ l'opr\'{e}%
rateur d\'{e}fini par $V(S)=Q_{1}\circ S\in \mathcal{L}(X,Y),$ $S\in 
\mathcal{L}(X,Y^{\ast \ast })$ et $Q:K(X,Y)^{\ast \ast }\rightarrow K(X,Y)$
l'op\'{e}rateur d\'{e}fini par $Q(R)=P\left[ V(\xi _{R})\right] ,$ $R\in
K(X,Y)^{\ast \ast }.$ On remarque que $Q$ est une projection continue sur $%
K(X,Y).\blacksquare $

D\'{e}monstration du th\'{e}or\`{e}me \ref{ert}.

D'apr\`{e}s \cite[Th4.1.1]{Tab}, il existe un espace de Banach $Z$ tel que $%
Z^{\ast }=\ell ^{1}$ et $K(Z)$ soit compl\'{e}ment\'{e} dans $\mathcal{L}%
(Z). $ Remarquons que l'espace $Z$ n'est pas r\'{e}flexif, comme $Z^{\ast
\ast }=\ell ^{\infty }$ est un espace de Grothendieck \cite{Groth} $Z$ ne
peut pas \^{e}tre compl\'{e}ment\'{e} dans son bidual. D'autre part, $%
Z^{\ast }=\ell ^{1}$ \ a la propri\'{e}t\'{e} de l'approximation born\'{e}e,
d'apr\`{e}s \cite[Chap.VIII-3]{Du} $Z$ a la propri\'{e}t\'{e} de
l'approximation born\'{e}e. En appliquant la proposition \ref{rutt}, on voit
que $K(Z)$ n'est pas compl\'{e}ment\'{e} dans son bidual.$\blacksquare $

\begin{problem}
\label{vc}Existe-il deux espaces de Banach $X$ et $Y$ tel que $K(X,Y)$ soit
compl\'{e}ment\'{e} dans son bidual et $\mathcal{L}(X,Y)\neq K(X,Y)?$
\end{problem}

\begin{remark}
\label{rrte}Si $Y$ ou $X^{\ast }$ a la propri\'{e}t\'{e} de l'approximation
born\'{e}e, il existe une projection de $\mathcal{L}(X,Y)^{\ast \ast }$ sur $%
K(X,Y)^{\ast \ast }.$
\end{remark}

Preuve. En effet,

Supposons que $Y$ a la propri\'{e}t\'{e} de l'approximation born\'{e}e. Soit 
$P_{1}:K(X,Y)^{\ast \ast \ast \ast }\rightarrow K(X,Y)^{\ast \ast }$ une
projection continue telle que $(P_{1}V,\xi )=(V,\xi )$ pour $\xi \in
K(X,Y)^{\ast }$ et $V\in K(X,Y)^{\ast \ast }.$

L'application $P_{1}\circ U^{\ast \ast }:\mathcal{L}(X,Y)^{\ast \ast
}\rightarrow (K(X,Y)^{\ast \ast }$ est une projection continue ($U$ est l'op%
\'{e}rateur d\'{e}fini auparavant).$\blacksquare $

\begin{remark}
\label{hhy}Soit $X$ un espace de Banach. Supposons qu'il existe une
projection continue $P:L^{\infty }(\mathbb{T},X)\rightarrow L^{\infty }(%
\mathbb{T})\overset{\vee }{\otimes }X.$ Alors $X$ est de dimension finie.
\end{remark}

Preuve. En effet,

Comme $L^{\infty }(\mathbb{T},X)$ contient une copie (isom\'{e}trique) compl%
\'{e}ment\'{e}e de $\ell ^{\infty }(X)$ et $L^{\infty }(\mathbb{T})\overset{%
\vee }{\otimes }X$ contient une copie compl\'{e}ment\'{e}e de $\ell ^{\infty
}\overset{\vee }{\otimes }X,$ alors $K(\ell ^{1},X)=\ell ^{\infty }\overset{%
\vee }{\otimes }X$ est compl\'{e}ment\'{e} de $\mathcal{L}(\ell ^{1},X)=\ell
^{\infty }(X).$ D'apr\`{e}s \cite[th.6]{Kalt}, $X$ est de dimension finie.$%
\blacksquare $

M.daher@orange.fr

\bigskip

\end{document}